\documentclass[12pt]{article}

      \usepackage{latexsym}
         \usepackage[reqno, namelimits, sumlimits]{amsmath} 
         \usepackage{amssymb, amsfonts}
         \usepackage{amsthm} 
\newtheorem{Theorem}{Theorem}[section]
\newtheorem{Lemma}{Lemma}[section]
\newtheorem{Proposition}{Proposition}[section]
\newtheorem{Remark}{Remark}[section]
\newcommand{\rth}{R^3}
\newcommand{\rtmiz}{R^2\times(-\infty,0)}
\newcommand{\rthmiz}{R^3\times(-\infty,0)}
\newcommand{\rnmiz}{R^n\times(-\i,0)}
\newcommand{\Om}{\Omega}
\newcommand{\Omt}{\Omega\times(0,T)}
\newcommand{\rnt}{R^n\times(0,T)}

\newcommand{\rtht}{R^3\times(0,T)}
\newcommand{\rthtp}{R^3\times(0,T')}
\newcommand{\vt}{\mathcal{V}_T}
\renewcommand{\div}{\operatorname{div}}
\newcommand{\curl}{\operatorname{curl}}
\newcommand{\heat}{(\partial_t-\Delta)}
\newcommand{\pd}[1]{{\frac{\partial}{\partial x_{#1}}}}
\newcommand{\kse}[1]{{#1}^{(k)}}
\newcommand{\irn}{\int_{R^n}}
\newcommand{\izt}{\int_0^t}
\newcommand{\izT}{\int_0^T}
\newcommand{\lirn}{L^{\i}(R^n)}

\newcommand{\vf}{\ensuremath{\varphi}}

\newcommand{\be}{\begin{equation}}
\newcommand{\ee}{\end{equation}}
\newcommand{\lxt}[2]{L_{x,\,t}^{#1}}
\newcommand{\lxts}[2]{L_{x,\,t}^{#1,\,#2}}
\renewcommand{\i}{\infty}
\newcommand{\liot}{\lxt\i\i(\rnt)}
\newcommand{\liotp}{\lxt\i\i(R^n\times(0,T'))}
\newcommand{\phie}{\phi_{\varepsilon}}
\newcommand{\ve}{\varepsilon}
\newcommand{\qr}{Q((\bar x,\bar t),R)}
\newcommand{\om}{\omega}
\newcommand{\g}{\gamma}

\renewcommand{\t}{\theta}
\newcommand{\q}[2]{{#1}_{#2}}
\newcommand{\dq}[3]{{#1}_{{#2},\,{#3}}}
\newcommand{\dtq}[3]{{#1}_{{#2}\,{#3}}}
\renewcommand{\o}{\omega}
\newcommand{\rut}{r\q u \t}
\newcommand{\otr}{\frac{\,\q \o \t}{r}}
\newcommand{\er}{\frac{\partial\,\,}{\partial r}}
\newcommand{\et}{\frac{\partial}{r\partial\t}}
\newcommand{\ez}{\frac{\partial\,\,}{\partial z}}
\newcommand{\fl}{f^{\lambda}}
\newcommand{\ul}{u^{\lambda}}

\renewcommand{\d}{\delta}
\newcommand{\ro}{\mathcal{R}_1}
\renewcommand{\r}{\sqrt{x_1^2+x_2^2}}
\newcommand{\bxt}{\bar x_3}
\newcommand{\itix}[1]{\int_{-\i}^0\int_{R^3}{#1}\,dx\,dt}
\newcommand{\vk}{v^{(k)}}
\newcommand{\wk}{w^{(k)}}
\newcommand{\ck}{\mathcal{C}_k}

\newenvironment{myproof}{\begin{proof}[{\rm\bf Proof.}]}{\end{proof}}

\newcommand{\R}{{\mathbb R}}

\title{
Liouville theorems for the Navier-Stokes equations and applications }

 \author{G. Koch\footnote{University of Chicago}\and 
  N. Nadirashvili\footnote{CNRS Laboratoire d'Analyse, Marseille} \and G. Seregin\footnote{Oxford University} 
\and V. \v Sver\'ak\footnote{University of Minnesota. Supported in part by NSF
  Grant DMS-0457061}} 

\begin{document}

\numberwithin{equation}{section}

\maketitle

\begin{abstract} We  study bounded ancient solutions of the Navier-Stokes equations. These are solutions with bounded 
velocity defined
in $R^n\times(-\i,0)$.  In two space dimensions we prove that such solutions are either constant or of the form
$u(x,t)=b(t)$, depending on the exact definition of admissible solutions. The general three dimensional problem
seems to be out of reach of existing techniques, but partial results can be obtained in the case of
axi-symmetric solutions. We apply these results to some scenarios of potential singularity formation for axi-symmetric
solutions, and obtain extensions of results in a recent paper by Chen, Strain, Tsai and Yau \cite{ChenStrainTsaiYau}. 
\end{abstract}

\section{Introduction}\label{sect1}
It is a well-known principle in the regularity theory of PDE that re-scaling procedures
are very useful in studying potential singularities. For example, for a minimal
surface $\Sigma\subset R^n$ for which $0\in\Sigma$ is a singular point, one
should look at the surfaces $\lambda\Sigma$ in the limit $\lambda\to\infty$,
see for example \cite{Giusti}. This ``blow-up" procedure, probably first introduced
by DeGiorgi in his study of minimal surfaces, has become indispensable 
in the study of singularities of various geometric equations
(see for example \cite{ Hamilton, SchoenUhlenbeck, Struwe}). Analogous
ideas were introduced in the study of many other classes of equations, such as semi-linear 
heat equations \cite{KohnGiga}, the Navier-Stokes equations 
\cite{ChenStrainTsaiYau, SereginSverak} and dispersive equations \cite{ KenigMerle, Struwe}, 
to name a few.
The blow-up procedure can be compared to infinite magnification and therefore
typically produces solutions of the original equation which are in some sense global.
The study of such global solutions is often a valuable stepping stone towards  understanding the structure
of potential singularities (or the absence of singularities).
In this paper we address some of these issues in the context of the Navier-Stokes equations
\be
\label{nse1}
\begin{array}{rcl}
u_t+u\nabla u +\nabla p -\Delta u & = & 0 \\
\div u & = & 0\,\,.
\end{array}
\ee
The scaling symmetry of the equations is $u(x,t)\to\lambda u(\lambda x,\lambda^2 t)$,
$p(x,t)\to \lambda^2p(\lambda x,\lambda^2 t)$ and can be used to ``zoom in" on a 
solution near a potential singularity. There are some free parameters in this process,
as we can choose where exactly (in space and time) we magnify (it does not have to be
exactly at a singularity, it can for example be just before the singularity occurs),
and which properties of the re-scaled solutions we wish to control. 
In this paper we study the situation in which we choose the $L^\i-$ norm of the re-scaled
velocity on a certain time interval as the parameter we wish to control.   
The pressure will play no explicit role in the process.
As we will see in Section~\ref{liouville}, this leads naturally to the following global problem:

\noindent
{\sl Characterize solutions of \eqref{nse1} in $R^n\times(-\i,0)$ with (globally) 
bounded velocity $u$.}

Following \cite{Hamilton}, we will call solutions defined in $R^n\times(-\i,0)$
{\it ancient solutions}. Stated in this terminology, we are interested in  ancient solutions
of \eqref{nse1} with bounded velocity.  A first guess might be that such solutions 
should be constant. To make this a plausible conjecture, one must be slightly more precise. 
Equation \eqref{nse1} has trivial non-constant solutions of the form $u(x,t)=b(t)$, $p(x,t)=-b'(t)x$
and so we need a definition of solutions which would eliminate these ``parasitic solutions".
The right definition seems to be that of a {\it mild solution} (see Section~\ref{linear}),
which was probably introduced in \cite{Kato}. (Implicitly it is already used in Leray's paper \cite{Leray}.) 
Another natural definition often used in the literature is that of a {\it weak solution}, 
also essentially introduced in Leray's paper \cite{Leray}, which 
is defined  using divergence-free test functions, see Section~\ref{linear}.
This notion of solution does allow the parasitic solutions above.
In these settings, the best possible result one can hope for which is consistent
with what is known about the equations would be that any ancient mild solution with bounded velocity is
constant and any ancient weak solution with bounded velocity is of the form $u(x,t)=b(t)$.
We will prove that this is indeed the case in dimension two and also in the case
of axi-symmetric fields in dimension three, if some additional conditions are satisfied
(see Section~\ref{liouville}). The case of general three-dimensional fields is, as far as we know,
completely open. In fact, it is open even in the steady-state case ($u$ independent of $t$).

The methods we use in the proofs of these results are elementary. The key component of the proof
in dimension two is the use of the vorticity equation: 
\be
\label{vorticity1}
\om_t+u\nabla \om=\Delta \om
\ee
This is a scalar equation and $\om$ satisfies the Harnack inequality (see e.\ g.\ \cite{Safonov}),
which can be used to show that if $\om\ne0$, then in large areas of space-time $\om$ has to be almost equal
to its maximum/minimum. (In fact, the strong maximum principle together with standard compactness results
is sufficient to prove this.)  This turns out to be incompatible with the boundedness of $u$.
(One might speculate that with the condition $\div u=0$,  a Liouville theorem might be true
for \eqref{vorticity1} at a linear level, without using the relation between $u$ and $\om$.
This, however, appears to be false -- see \cite{UnpublishedCalculation}.) 
The ideas behind the proofs of the results for axi-symmetric fields in dimension three are similar.
In each case there is a scalar quantity satisfying a maximum principle which is  used in a way
similar to the two-dimensional case. The quantities we use and the corresponding
maximum principles are all classical.

There is a technical component in the  proofs, since one needs to establish that the solutions
we work with have sufficient regularity. This part is more or less standard,  and we use elementary techniques
based on explicit representation formulae to establish the required properties. 

In the last section we use the Liouville theorems of Section~\ref{liouville} to obtain results limiting
the types of singularities which may occur in axi-symmetric  solutions of the Navier-Stokes equations. These results
are inspired by the recent paper \cite{ChenStrainTsaiYau}, where significant progress
in the study of the axi-symmetric case was made using methods quite different from the ones presented
here. Our results on axi-symmetric singularities 
address some questions which were left open in \cite{ChenStrainTsaiYau}.
Very recently we learned that the authors of \cite{ChenStrainTsaiYau} have independently proved
results similar to those in Section~\ref{axisymsing} using their own methods.
Their paper \cite{ChenStrainTsaiYau2} on the subject is expected to appear soon.

 It is  known
that axi-symmetric solutions with no swirl have to be regular, see \cite{Ladyzhenskaya, Yudovich}. (We recall that the ``no swirl" condition
means that in cylindrical coordinates $(r,\t,z)$ -- see \eqref{cylindrical} -- the $\q u\t-$component
of the velocity vanishes.) However, the case of non-zero swirl is open at the time of this writing.
We will prove that, under natural assumptions, every potential singularity of an axi-symmetric solution
has to be of type II, in the sense of \cite{Hamilton}.
We recall that a singularity of a Navier-Stokes solution $u$ at time $T$ is called type I if
$$
\sup_x|u(x,t)|\le \frac{C}{\sqrt{T-t}}\,\,
$$
for some $C>0$. By definition, a type II singularity is any singularity which is not of type I.
A blow up of $u$ by a type II singularity is sometimes called {\it slow blow-up}, see e.\ g.\
\cite{Hamilton}.  Therefore we can rephrase our result by saying that if an axi-symmetric solution 
develops a singularity, it can only be through slow blow-up.
We remark that Leray proved in \cite{Leray} that if $u$ develops a singularity at $T$,
then 
$$
\sup_x|u(x,t)|\ge\frac{\ve_1}{\sqrt{T-t}}
$$
for some $\ve_1>0$. Also, the rate $\frac1{\sqrt{T-t}}$ would be the blow-up rate 
of a self-similar singularity. (It is known that these do not exist, see \cite{NRS, Tsai}.)

It is worth mentioning that although our results are obtained by methods which are
more or less elementary, it seems that some of them are out of reach of the usual methods
 used in the theory of the Navier-Stokes equations, such as energy methods or perturbation
 analyses in various function spaces. This is because some special properties of solutions of
 scalar equations, although simple, cannot be detected at the broad level 
 at which the usual methods used for Navier-Stokes are applied. 
 A similar situation appears in the proof that Leray's self-similar singularities do not exist,
 see \cite{NRS, Tsai}, where a (non-classical) scalar quantity satisfying an elliptic
 equation is used.
 At the time of this writing, there is no known similar quantity for the general three-dimensional problem.

\noindent

\section{Preliminaries}\label{preliminaries}
Let $\Om\subset R^n$ be a bounded domain and let $T>0$.
We consider the parabolic equation in $\Omt$  of the form
\be
\label{parequation}
u_t+a(x,t)\nabla u-\Delta u=0\,\,,
\ee
with $a\in\lxt\i\i(\Omt)$. A suitable notion of a solution is for example
a weak solution. By definition, $u$ is a weak solution of \eqref{parequation}
if $u$ and $\nabla_x u$ (the distributional derivative) belong
to $({\lxt 2 2})_{\rm loc}(\Omt)$ and the equation is satisfied in distributions.
It then follows from standard regularity that  in fact $u_t$ and $\nabla_x^2u$ belong to
$({\lxt p p})_{\rm loc}(\Omt)$ for every $p\in(0,\i)$ and the equation is satisfied
pointwise  almost everywhere in $\Omt$. See for example \cite{LSU}.
Therefore there is no difference between weak solutions and strong solutions,
and we can just use the term ``solution" in the context of \eqref{parequation}.
We recall that the ``parabolic boundary" of $\Omt$ is 
$\partial_{\rm par}(\Omt)=(\bar\Omega\times\{0\})\cup (\partial\Om\times[0,T])$.
When $x\in\Om$, the space-time points $(x,T)$ belong to the
``parabolic interior" of $\Omt$ and $u(x,T)$ is well-defined. 
We recall that the solutions of \eqref{parequation} satisfy the
{\it strong maximum principle}: If $u$ is a bounded solution in $\Omt$
such that $u(\bar x,T)=\sup_{\Omt} u$ for some $\bar x\in\Om$, then $u$ is constant
in $\Omt$. In fact, a much stronger statement is true: non-negative solutions of 
\eqref{parequation} satisfy the parabolic Harnack inequality, see
for example \cite{Safonov}. The Harnack inequality immediately implies 
the strong maximum principle.  For our purposes in this paper the strong
maximum principle is sufficient -- we will not need the full strength
of the Harnack inequality. Our key tool will be the following lemma which
essentially says that the statement of the strong maximum principle
is in some sense stable under perturbations. (This stability can be made much more 
precise with the Harnack inequality.) The lemma is certainly known in one form 
or another, but we were unable to locate in the literature the precise statement we
need.

\begin{Lemma} 
\label{weakharnack}
Let us consider  equation \eqref{parequation} with bounded 
measurable
coefficient $a$ in $\Om\times(0,T)$.
Let $K$ be a compact subset of $\Om$,  $\Om'\subset\bar\Om'\subset\Om$
and  $\tau>0$. Then for each $\ve>0$ there exists $\delta>0$,
$\delta=\delta(\Om,\Om',K,T,||a||_{\lxt\i\i},\tau,\ve)$ such that if $u$ is a bounded solution
of \eqref{parequation} with $\sup_{\Omt}|u|= M$ and
$\sup_{x\in K}u(x,T)\ge M(1-\delta)$, then
$u(x,t)\ge M(1-\ve)$ in $\Om'\times(\tau,T)$.
\end{Lemma}

\begin{myproof} We can take $M=1$ without loss of generality.
Assuming the statement fails for some $\ve>0$, there must exist a sequence of 
coefficients $a^{(k)}$,
solutions $u^{(k)}$ of \eqref{parequation} with $a=a^{(k)}$, and  points $x_k\in K$ 
and $(y_k,t_k)\in\Om'\times(\tau,T)$ such that $|a^{(k)}|\le C$, 
$|\kse u|\le 1$, $\kse u(x_k,T)\to 1$ and  $\kse u(y_k,t_k)\le 1-\ve$.
We can assume, after passing to a subsequence, that
$\kse a$ converge  weakly$^*$ in $\lxt\i\i$ to $\bar a$,
$\kse u$ converge locally uniformly  in 
$\Omt$ to $\bar u$, $x_k\to\bar x\in K$ and $(y_k,t_k)\to (\bar y,\bar t)\in
\bar\Om'\times[\tau, T]$. The regularity properties of solutions
of $\eqref{parequation}$ discussed above imply that $\bar u$ solves
\eqref{parequation} with $a=\bar a$, $|\bar u|\le 1$ in $\Omt$,
$\bar u(\bar x,T)=1$ and $\bar u(\bar y,\bar t)\le 1-\ve$. 
This, however, is impossible due to the strong maximum principle.
\end{myproof}

\section{Bounded solutions of the linear Stokes \\ problem}\label{linear}
Let us first recall some basic facts about the Cauchy problem for the linear
Stokes system, with $u=(u_1,\dots,u_n)\colon R^n\times(0,\infty)\to R^n$
and the right-hand side in divergence form:
\begin{eqnarray}
%\label{linstokes}
\left.
\begin{array}{rcl}
u_t+\nabla p -\Delta u & = & \pd k f_k \label{linstokeseq}\\
\div u & = & 0
\end{array}
\right\}&  & 
\mbox{in $R^n\times(0,\infty)$}\\
\begin{array}{ccl}
u(\,\cdot\,,0) & = & u_0\,\,\,\,\,\,\,\,\,\,\,\,
\end{array} & & \, \mbox{in $R^n$} \label{ic}
\end{eqnarray}

Here $f_k=(f_{1k},\dots,f_{nk})$ for $ k=1,\dots,n$. 
Denoting by $P$ the Helmholtz projection of vector fields on div-free
fields and by $S$ the solution operator of the heat equation,
we have the well-known representation formula
\begin{equation}
\label{abstractrep}
u(t)=S(t)u_0+\izt S(t-s)P\pd k f_k(s)\,ds\, ,
\end{equation} 
where, as usual, $u(t)$ denotes the function $u(\,\cdot\,,t)$, etc.

This can be written more concretely in terms of the kernel
$$
K_{ij}(x,t)=(-\delta_{ij}\Delta+\frac{\partial^2}{\partial x_i\partial x_j})\Phi(x,t),
$$
where the ``generating function" $\Phi$ is defined in terms of the fundamental 
solution of the Laplace operator
$G$ and the heat kernel $\Gamma$ by
\be
\label{generatingfunction}
\Phi(x,t)=\irn G(y)\Gamma(x-y,t)\,dy,
\ee
which is the same as
$$
\Phi(\,\cdot\,,t)=S(t)G.
$$
See for example \cite{Oseen}.
Letting
$$
K_{ijk}=\pd k K_{ij}\,,
$$
we can re-write \eqref{abstractrep} as
\begin{equation}
\label{rep}
u_i(x,t)=\irn \Gamma(x-y,t)u_{0i}(y)+\izt\irn K_{ijk}(x-y,t-s)f_{jk}(y,s)\,dy\,ds\,.
\end{equation}
Note also the obvious estimates
\begin{equation}
|K_{ij}(x,t)|\le \frac{C}{(|x|^2+t)^{\frac{n}{2}}}
\end{equation}
and
\begin{equation}
\label{decay}
|K_{ijk}(x,t)|\le \frac{C}{(|x|^2+t)^{\frac{n+1}{2}}}\,\,.
\end{equation}
As a consequence of \eqref{decay}, the expression \eqref{rep} is well-defined for
$f\in\lxt\i\i$. We remark that, in contrast, solutions of 
\begin{eqnarray}
\label{linstokesnd}
\left.
\begin{array}{rcl}
u_t+\nabla p -\Delta u & = &  f \\
\div u & = & 0
\end{array}
\right\}&  & 
\mbox{in $R^n\times(0,\infty)$}\\
\begin{array}{ccl}
u(\,\cdot\,,0) & = & u_0\,\,\,
\end{array} & & \, \mbox{in $R^n$}
\end{eqnarray}
are not well defined for $f\in\lxt\i\i$, although the ambiguity is small. This can also be seen
without using the explicit form of the kernel, in the following way:
One can write, for each $t$, the Helmholtz decomposition of $f(x,t)$ as
$f(x,t)=Pf(x,t)+\nabla_x\phi(x,t)$. The projection $P$ can be naturally defined on $L^{\i}(R^n)$
(which is mapped by $P$ into  $BMO(R^n)$ ) only modulo  constants, which creates an ambiguity. However,
if the right-hand side is in divergence form, this ambiguity is cancelled by the extra
derivative.

By definition, a {\it mild solution} of the Cauchy problem \eqref{linstokeseq} and \eqref{ic} is a function $u$ defined by 
the formula \eqref{rep}. We note that this definition does not involve the pressure.
One can obtain (formally) an explicit formula for the pressure, but, unlike the formula for the velocity field $u$,
  it defines $p$ only modulo a function of $t$ (constant in $x$ for 
 each $t$) when $f_k$ is in $\lxt\i\i$. 
 
 The definition of mild solutions immediately implies their uniqueness. 
 Also, we have standard estimates for $u$ in terms of $f=(f_1,\dots,f_n)=(f_{ij})_{i,j=1}^n$.
 In particular, for $u_0=0$ we have the  estimates
 \begin{eqnarray}
 %\label{estimates}
 ||u||_{C^{\alpha} _{\rm par}( Q(z_0,R) )} & \le &  C(\alpha,R)||f||_{\lxt\i\i(R^n\times(0,T))} \quad \mbox{and}\label{estimates1}\\
 ||\nabla_x u||_{\lxt p p(Q(z_0,R))} & \le & C(p,R) ||f||_{\lxt\i\i(R^n\times(0,T))}\,\,\label{estimates2}
 \end{eqnarray}
 for any $\alpha\in(0,1)$ and $p\in(1,\i)$, where $Q(z_0,R)=Q((x_0,t_0),R)=B(x_0,R)\times (t_0-R^2,t_0)$ is 
 any parabolic ball contained in $R^n\times(0,T)$.  The space $C^{\alpha}_{\rm par}$ is defined by means of the parabolic
  distance $\sqrt{|x-x'|^2+|t-t'|}$.
  
  Taking difference quotients, we see that, for $u_0=0$, we have similar estimates for spatial derivatives:
  \begin{eqnarray}
  \label{estimatesder}
 ||\nabla_x ^k u||_{C^{\alpha}_{\rm par}(Q(z_0,R))} & \le &  C(\alpha,R)||\nabla_x^kf||_{\lxt\i\i(R^n\times(0,T))} \quad \mbox{and}\\
 ||\nabla_x^{k+1} u||_{\lxt p p(Q(z_0,R))} & \le & C(p,R) ||\nabla_x^k f||_{\lxt\i\i(R^n\times(0,T))}\,\,
 \end{eqnarray}
 Moreover, a routine inspection of representation formula \eqref{rep} shows that, when $u_0=0$, the time derivative 
 satisfies, for $k=0,1,\dots$,
 \be
 \label{estimatetime}
 ||\nabla_x^{k} u_t||_{\lxt\i\i(R^n\times(0,T))}\le C(T,k)||\nabla_x^{k+2} f||_{\lxt\i\i(R^n\times(0,T))}.
 \ee
 We sketch the calculation leading to the last estimate in the case $k=0$ for the convenience of the reader:
 Clearly it is enough to estimate $|u_t(0,t)|$. Let $\Phi$ be the generating function defined in
 \eqref{generatingfunction}, which will be considered as a function of $R^n\times R$, with $\Phi=0$
 for negative values of $t$. We can write
 \be
 \label{convolutionrep}
 u_i=(L_{ijk}\Phi)*f_{jk}, 
 \ee
 where $L_{ijk}$ is a homogeneous constant coefficient operator in $x$ of order 3 and $*$ denotes space-time 
 convolution. 
 Applying the heat operator to \eqref{convolutionrep} we can write, with a slight abuse of notation,
 \be
 \label{heatu}
 \heat u_i=(L_{ijk}\heat\Phi)*f_{jk}=(L_{ijk}G(x)\delta(t))*f_{jk},
 \ee
 where $G$ is the fundamental solution of the Laplacian and $\delta(t)$ is the
 Dirac distribution in $t$.
 We consider a smooth cut-off function $\eta=\eta(x)$ on $R^n$
 with $\eta=1$ in the unit ball $B(0,1)$ and $\eta=0$ outside of $B(0,2)$ and set
 $f'=\eta f$, $f''=(1-\eta) f$.
 Let us first look at $u'_i$, the contribution to $u_i$ in 
 \eqref{heatu} coming from $f'$. We can move two derivatives from $L_{ijk}$ to $f'_{jk}$
 to obtain an estimate of $\heat u_i'(0,t)$ in terms of the $\lxt \i\i-$norm of the second derivatives of $f'_{jk}$.
 The estimate of $\heat u_i''(0,t)$ (with the obvious meaning of $u_i''$) is even simpler, since 
 $L_{ijk}G$ is integrable in $R^n\setminus B(0,1)$ and therefore $\heat u_i''(0,t)$ can be estimated
 in terms of the $\lxt\i\i-$norm of $f''_{jk}$. 
 Once we have the estimate for $\heat u$, the estimate for $u_t$ follows
 from \eqref{estimatesder}.

 To define the notion of a weak solution of equation \eqref{linstokeseq}, we follow the standard procedures and
 introduce the space $\vt$ of smooth compactly supported div-free vector fields $\vf\colon\rnt\to R^n$. 
 We then say that a bounded measurable vector field $u\colon\rnt\to R^n$ is a {\it weak solution}
 of \eqref{linstokeseq} if $\div u = 0$ in $\rnt$ (in the sense of distributions) and
 $\izT\irn u(\varphi_t+\Delta\varphi)\,dx\,dt=\izT\irn f_k\pd k\vf \,dx\,dt$ for each $\varphi\in\vt$. 
 
 \begin{Lemma}
 \label{decomposition}
 For a fixed $f\in\lxt\i\i(\rnt)$ let $u\in\lxt\i\i(\rnt)$ be any weak solution of \eqref{linstokeseq} in $\rnt$, 
 and denote by $v$  the mild solution of the Cauchy problem \eqref{linstokeseq} and \eqref{ic} with
 $u_0=0$. Then $u(x,t)=v(x,t)+w(x,t)+b(t)$, where $w$ satisfies the heat equation $w_t-\Delta w = 0$ in $\rnt$ and
 $b$ is a bounded measurable $R^n-$valued function on $(0,T)$. Moreover, we have the estimates
 \begin{eqnarray}
 \label{estimatesbw}
 ||w||_{\liot} & \le & C(T) ||u||_{\liot} \quad \mbox{and} \\
 ||b||_{L^\i(0,T)} & \le & C(T) ||u||_{\liot}\,.
 \end{eqnarray}
 \end{Lemma}
 \begin{myproof}
 In view of estimates \eqref{estimates1} it is enough to consider only the case $f=0$. Let $\phi\colon R^n\times R\to R$
 be a mollifyer compactly supported in $R^n\times(-1,0)$, 
 $\phie(x,t)=\ve^{-(n+1)}\phi(x/\ve,t/\ve)$, and let $u_{\ve}\colon R^n\times(0,T-\ve)$ be defined by
$ u_\ve=\phie*u$ (space-time convolution). Let $w_\ve$ be the solution
 of the heat equation in $\rnt$ with  initial datum $w_\ve(x,0)=u_\ve(x,0)$.
 The (smooth and bounded) function $h_\ve=\curl(u_\ve-w_\ve)$ satisfies the heat equation in
 $R^n\times(0,T-\ve)$ with initial datum $h_\ve(x,0)=0$ and therefore it must vanish.
 Since bounded solutions of the system $\curl z=0$ and $\div z=0$ in $R^n$ are constant
 by Liouville's theorem, we see that $u_\ve(x,t)-w_\ve(x,t)=b_\ve(t)$ for a suitable
 $b_\ve\colon(0,T-\ve)\to R^n$. By compactness properties of families of bounded solutions of the heat
 equation we see that if  $\ve\to 0$ along a suitable sequence, the functions $b_\ve$ converge a.\ e.\ 
 to an $L^\infty$ function
 $b\colon(0,T)\to R^n$. The estimates follow from the constructions.
 \end{myproof}
 \begin{Remark}
 In the above decomposition, the function $v$ is of course uniquely determined by $f$, whereas
 the functions $w$ and $b$ are determined up to a constant (independent of time). In other words,
 the (distributional) derivative $b'(t)$ is uniquely determined by $u$ and $f$.
 \end{Remark}

\section{Bounded solutions of Navier-Stokes }\label{nsreg}
Let us now consider the Cauchy problem for the Navier-Stokes equations:
\begin{eqnarray}
%\label{NS}
\left.
\begin{array}{rcl}
u_t+u\nabla u +\nabla p -\Delta u & = & 0 \\
\div u & = & 0
\end{array}
\right\}&  & 
\mbox{in $R^n\times(0,\i)$} \label{NSeq}\\
\begin{array}{ccl}
u(\,\cdot\,,0) & = & u_0\,\,
\end{array} & & \quad \mbox{in $R^n$}\label{NSic}
\end{eqnarray}
The considerations of the previous section can be repeated with $f_k=-u_ku$. In particular,
a function $u\in\lxt\i\i(\rnt)$ is defined to be (i) a {\it mild solution} of the Cauchy problem
\eqref{NSeq} and \eqref{NSic} if \eqref{rep} is valid with $f_k=-u_ku$ and (ii) a
{\it weak solution} of equation \eqref{NSeq} in $\rnt$ if   
$\div u = 0$ in $\rnt$ (in the sense of distributions) and
 $\izT\irn u(\varphi_t+\Delta\varphi)\,dx\,dt=\izT\irn -u_ku\pd k\vf \,dx\,dt$ for each $\varphi\in\vt$. 
 
 \begin{Remark} It is obvious that the notions of weak solution and mild solution 
 are also well defined under 
 the assumption that
 $u\in\lxt\i\i(R^n\times(0,T')$ for each $T'<T$ (with the possibility that 
 $||u||_{\lxt\i\i(R^n\times(0,T'))}\to\i$ as $T'\nearrow T$). This is a natural
 setting in which potential singularities of solutions of the Cauchy problem can be studied.
 Even if one considers the Cauchy problem for $u_0$ in spaces other
 than $\lirn$, such as $L^n(R^n)$ (\cite{Kato}) or $BMO^{-1}(R^n)$ (\cite{KochTataru}),
 the local-in-time solution $u\colon\rnt\to R^n$ which is constructed for
 $u_0$ in these spaces typically belongs to $\lxt\i\i(R^n\times(\tau,T-\tau))$ for
 any $\tau>0$. 
 \end{Remark}

 The existence and uniqueness of {\it local-in-time} mild solutions 
 of the Cauchy problem \eqref{NSeq} and \eqref{NSic} with $u_0\in L^{\infty}$
 was  addressed in \cite{Giga}.  We briefly outline a slightly modified  
 approach using standard perturbation theory. We define the bilinear
 form $B\colon\liot\times\liot\to\liot$ by 
 \be
 \label{blform}
 B(u,v)_i(x,t)=\izt\irn -K_{ijk}(x-y,t-s)u_k(y,s)v_j(y,s)\,dy\,ds\,\,,
 \ee
 and we  denote by $U$ the heat extension of the initial datum $u_0$.
 The equation for $u$ then becomes 
 \be
 \label{abstracteq}
 u=U+B(u,u)
 \ee
 and can be solved in $\liot$ for sufficiently small $T$ by a fixed point argument, 
 since  estimate \eqref{decay} easily implies
 \be
 \label{Bestimate}
 ||B(u,v)||_{\liot}\le C\sqrt{T}||u||_{\liot}||v||_{\liot}.
 \ee
 We remark that \eqref{estimates1} implies that the solutions of \eqref{abstracteq} have enough regularity 
 to allow us to treat
 \eqref{abstracteq} as an ODE in $t$, without making  assumptions about $u$ other than $u\in\liot$.
 
 We recall now the regularity properties of mild solutions in $\liot$. The following (optimal) result
 will not be needed here in its full generality, but we feel it is still worth mentioning:
\begin{Proposition}\label{smoothing}
Let $u\in\liot$ be a mild solution of \eqref{NSeq} and \eqref{NSic} with $u_0\in L^{\infty}$.
Then for $k,l=0,1,\dots$ the functions
$t^{\frac{k}{2}+l}\nabla_x^k\partial_t^l u$ are bounded and, for  $T'=\ve(k,l){||u_0||^{-2}_{L^\i(R^n)}}$ 
(where $\ve(k,l)>0$ is a small constant), we have
\be
\label{smoothnessest}
||t^{{\frac{k}{2}}+l}\nabla_x^k\partial_t^l u||_{\liotp}\le C(k,l)||u_0||_{L^\i(R^n)}\,\,.
\ee
\end{Proposition}
\begin{myproof}
This can be proved in the same way as the corresponding results in
\cite{GigaSawada}, \cite{DongDu} and \cite{GermainPavlovicStaffilani}, where the authors
work in function spaces other than $\lxt\i\i$. The key is an estimate of $B$ with the same form
as \eqref{Bestimate} but in spaces with norms given by the expression on the
left-hand side of \eqref{smoothnessest}. In the context of the $\liot-$based norms we use here,
the proof is in fact much simpler than in that of the spaces used in the above papers, 
due to the elementary nature of estimate \eqref{Bestimate}.
\end{myproof} 

\begin{Remark} Estimate \eqref{smoothnessest} says that the local-in-time 
smoothing properties of Navier-Stokes for $u_0\in L^\i$ are the same as those of the
heat equation. Since the solution $u$ is constructed essentially as a power series 
perturbation around the heat extension $U$ of $u_0$, this may not be surprising.
\end{Remark}

\begin{Lemma}
Let $u^{(k)}\in\liot$ be a sequence of mild solutions of $\eqref{NSeq}$ and \eqref{NSic} 
with initial conditions $u^{(k)}_0$. Assume 
$||u^{(k)}||_{\liot}\le C$ with $C$ independent of $k$. Then a subsequence
of the sequence $u^{(k)}$ converges locally uniformly in $R^n\times(0,T)$
to a mild solution $u\in\liot$ with initial datum $u(x,0)$ given by 
the weak$^*$ limit of a suitable subsequence of the sequence $u^{(k)}_0$.
\end{Lemma}

\begin{myproof} This is a routine consequence of \eqref{smoothnessest},
and the decay estimate \eqref{decay} for the kernel $K_{ijk}$. 
\end{myproof}

We now turn to regularity properties of bounded weak solutions.
Let $u\in\liot$ be a weak solution of \eqref{NSeq} in $\rnt$, and let
$M=||u||_{\liot}$.
Let $v$ be the mild solution of the linear Cauchy problem \eqref{linstokeseq}
and \eqref{ic} with $f_k=-u_k u$ and $u_0=0$. By Lemma~\ref{decomposition}
we can write $u=v+w+b$ with the $L^\i-$norms of $v,\,\,w$ and $b$ bounded
by $N=C_1(T)M^2+C_2(T)M$,  $w_t-\Delta w = 0$ and $b$ is a function
of $t$ only. Hence for $k=0,1,2,\dots$ and $\delta>0$ the derivatives
$\nabla_x^k(w+b)$ are bounded by $C(k,\delta)N$ in $R^n\times(\delta,T)$ by
Proposition~\ref{smoothing}. Moreover, we have the
$L^p-$estimate  \eqref{estimates2} for $\nabla_x v$.
Therefore $\omega=\curl u$ belongs to $\lxt p p(Q(z_0,R))$ for any
$p\in(1,\i)$ and any $Q(z_0,R)\subset R^n\times(\delta, T)$, with
\be
\label{omestimate}
||\omega||_{\lxt p p(Q(z_0,R))}\le C(p,\delta,R,M).
\ee

Following \cite{Serrin}, we can now use the equation for $\omega$
to obtain estimates for higher derivatives $\nabla_x^k u$. For $n=3$
the equation for $\omega$ is
\be
\label{omegaeq}
\omega_{i\,t}-\Delta\omega_i=\pd j(\omega_ju_i-\omega_i u_j)
\ee
and it is easy to check that in our situation this equation is satisfied
in the sense of distributions. Equation \eqref{omegaeq} gains $\omega$ one
spatial derivative  in $\lxt p p$. The standard bootstrapping
arguments and regularity estimates for harmonic functions now give
\be
\label{dxlpestimates}
 ||\nabla_x^k u||_{\lxt p p(Q(z_0,R))}\le C(k,\delta,R,M) 
\ee
for each $Q_(z_0,R)\subset R^n\times(\delta, T)$.
Therefore, using standard imbeddings, we have for $k=0,1,2\dots$
\be
\label{dxliestimates}
||\nabla_x^k u||_{\lxt\i\i(R^n\times(\delta,T))}\le C(k,\delta,T,M).
\ee
Finally, using \eqref{estimatetime} we also obtain for $k=0,1,2\dots$
\be
\label{dtestimate}
||\nabla_x^k \partial_t(u-b)||_{\lxt \i \i(R^n\times(\delta,T))}\le C(k,\delta,R,M).  
\ee
(We adopt the usual convention that the value of $C$ can change from line to line.)

\section{Liouville theorems}\label{liouville}
Let us first consider the Navier-Stokes equations in two space dimensions.
\begin{Theorem}
\label{nstd}
Let $u$ be a bounded weak solution of the Navier-Stokes equations
in $R^2\times(-\i,0)$. Then $u(x,t)=b(t)$ for a suitable bounded
measurable $b\colon(-\i,0)\to R^2$.
\end{Theorem}
\begin{myproof}
In two space dimensions the vorticity is a scalar quantity defined by
\be
\label{vortdef}
\omega=u_{2,1}-u_{1,2}\,\,,
\ee
where the indices after comma mean derivatives, i.\ e.\ $u_{2,1}=\pd 1 u_2$, etc.
By the results of Section~\ref{nsreg}, the function $\omega$ is uniformly bounded
together with its spatial derivatives. Moreover, its time derivative is also uniformly bounded. 
 The vorticity equation in dimension two is
\be
\label{vorteq}
\omega_t+u\nabla\omega-\Delta \omega = 0\,.
\ee
Let $M_1=\sup_{\rtmiz} \omega$, $M_2=\inf_{\rtmiz}\omega$ and assume that $M_1>0$.
By Lemma~\ref{weakharnack} there exist arbitrarily large balls
$Q_R=Q((\bar x,\bar t), R)=B(\bar x,R)\times(\bar t-R^2,\bar t)$ such that $\omega\ge M_1/2$ in 
$\qr$. For such balls we have 
\be
\label{volume}
\int_{Q_R} \omega\,dx\,dt\ge \pi M_1R^4.
\ee
On the other hand, denoting by $n$ the normal to the boundary of $B(\bar x, R)$, we can also write
%\be
%\label{surface}
%\int_{Q_R}\omega \,dx\,dt  =  \int_{Q_R}(u_{2,1}-u_{1,2})\,dx\,dt  =  
%  \int_{\partial B(\bar x,R)\times(\bar t-R^2,\bar t)}(u_2n_1-u_1n_2)\,ds\,dt\le CR^3\,.
%\ee
\be
\label{surface}
\begin{array}
{cl}
\int_{Q_R}\omega \,dx\,dt  & = \int_{Q_R}(u_{2,1}-u_{1,2})\,dx\,dt  =  \\   
& \int_{\partial B(\bar x,R)\times(\bar t-R^2,\bar t)}(u_2n_1-u_1n_2)\,ds\,dt\le  CR^3\,.
\end{array}
\ee
Clearly \eqref{volume} is not compatible with \eqref{surface}, unless $M_1\le0$.
In the same way we conclude that $M_2\ge 0$ and therefore $\omega$ must vanish identically.
Hence $\curl u = 0$ in $\rtmiz$ which, together with $\div u = 0$ and the boundedness of $u$,
implies (by the classical Liouville theorem for harmonic functions) that $u$ is constant in $x$ for each $t$.
\end{myproof}

It is not known if a result similar to Theorem~\ref{nstd} remains true in three spatial dimensions. In fact,
the problem is open even in the steady-state case. However, under the additional assumption that the
solutions are axi-symmetric, one can obtain some results which seem to be of interest.
We recall that a vector field u in $R^3$ is axi-symmetric if it is invariant under rotations about
a suitable axis, which is often identified  with the $x_3-$ coordinate axis.
In other words, a field $u$ is axi-symmetric if $u(Rx)=Ru(x)$ for every rotation $R$ of the form   
\begin{equation*}
R=\left(\begin{array}{ccc}\cos\alpha & -\sin\alpha\,\,\, & 0 \\ \sin\alpha &\, \cos\alpha & 0 \\
0 & 0 & 1\end{array}\right)\,\,.
\end{equation*}

In cylindrical coordinates $(r,\t,z)$ given by 
\be
\label{cylindrical}
x_1=r\cos\t,\quad x_2=r\sin\t,\quad x_3=z,
\ee
 the axi-symmetric
fields are given by $u=\q u r\er + \q u \t\et + \q u z \ez$, where the coordinate functions
$\q u r, \q u \t, \q u z$ depend only on $r$ and $z$. In these coordinates, the Navier-Stokes
equations become
\begin{eqnarray}
\dtq u r t + \q u r\dq u r r + \q u z\dq u r z -\frac{{\q {u} \t}^2}{r} + p_{,r} & = & \Delta \q u r - \frac{\q u r}{r^2} \label{ce1}\\
\dtq u \t t +\q u r\dq u \t r + \q u z\dq u \t z +\frac{\q u r\q u \t}{r} & = & \Delta \q u \t - \frac{\q u \t}{r^2}\label{ce2}\\
\dtq u z t + \q u r\dq u z r + \q u z\dq u z z + p_{,z} & = & \Delta \q u z\label{ce3}\\
\frac {(r\q u r)_{,r}}{r}+\dq u z z & = & 0\,\,,\label{c4}
\end{eqnarray}
where $\Delta$ is the scalar Laplacian (expressed in the coordinates $(r,\t,z)\,$), 
$\dq u r z$  denotes the partial derivative $\ez \q u r$, etc. 
The equation for $\q u \t$ is of special interest, as it is decoupled from the pressure.
The role of the non-linear terms in this equation can be seen by considering the inviscid
case (Euler's equations), wherein equation \eqref{ce2} is replaced by
\be
\dtq u \t t +\q u r\dq u \t r + \q u z\dq u \t z +\frac{\q u r\q u \t}{r}  = 0\,\,,
\ee
which is the same as
\be
\label{inviscid}
(\rut)_{\,t}+\q u r (\rut)_{,\,r} + \q u z (\rut)_{,\,z}=0.
\ee
Equation~\eqref{inviscid} says that the quantity $\rut$ ``moves with the flow".
This is a special case of Kelvin's law that the integral
of $u_idx_i$ along curves moving with the flow is constant.
In the situation considered here, the curves are circles centered at
the $x_3-$axis and lying in planes perpendicular to it.

In view of \eqref{inviscid}, it is natural to re-write
\eqref{ce2} as an equation for $\rut$:
\be
\label{ruteq}
(\rut)_{,t}+\q u r(\rut)_{,\,r} + \q u z (\rut)_{,\,z}=\Delta(\rut)-\frac{2}{r}(\rut)_{,\,r}
\ee

The infinitesimal version of Kelvin's law, which is  Helmholtz's law that vorticity
``moves with the flow" (for inviscid flows), gives in the case of axi-symmetric
flows without swirl ($u_\theta=0$) another quantity which moves with the flow, namely $\otr$. Here
$\om = \curl u$, as usual,  and in cylindrical
coordinates we write $\om = \q \om r\er + \q \om \t\et + \q \om z\ez$.
(For axi-symmetric flows without swirl we have $\q \om r=0$,
$\q \om z = 0$, and we can write $\om = \q\om\t\et$. Therefore the situation
is similar to two-dimensional flows.)

Hence for axi-symmetric solutions of Euler's equations without swirl we have
\be
\label{otequationeuler}
(\otr)_t+\q u r (\otr)_{,\,r}+\q u z (\otr)_{,\,z}=0\,.
\ee
This is nothing but the $\t-$component of the equation for $\om$, and can be
of course obtained by simple calculation, without any consideration of 
the Helmholtz law. 
For axi-symmetric solutions of Navier-Stokes without swirl the last equation becomes
\be
\label{otequationns}
(\otr)_t+\q u r (\otr)_{,\,r}+\q u z (\otr)_{,\,z}=\Delta(\otr)+\frac 2 r (\otr)_{,\,r}\,\,.
\ee

\begin{Remark} For a smooth vector field $u$, the apparent singularity of 
$\otr$ is only an artifact of the co-ordinate choice. The quantity 
$\otr$ is actually a smooth function, even across the $x_3-$axis,
as long as $u$ is smooth.
\end{Remark}

The diffusion term on the right-hand side of equation \eqref{otequationns}
can be interpreted as the $5-$dimensional Laplacian acting on 
$SO(4)-$invariant functions in $R^5$. We write $r=\sqrt{y_1^2+y_2^2+y_3^2+y_4^2}$,
$y_5=z$ and we note that for $\tilde f(y_1,\dots,y_5)=f(r,z)$ we have
\be
\Delta_y\tilde f(y_1,\dots,y_5)=
(\frac{\partial^2 f}{\partial r^2}+\frac{3\partial f}{r\partial r}+\frac{\partial ^2 f}{\partial z^2})(r,z)\,\,.
\ee
Therefore, with a slight abuse of notation,   we can write the equation \eqref{otequationns}
as
\be
\label{ote5}
(\otr)_t+\q u r (\otr)_{,\,r}+\q u z (\otr)_{,\,z}=\Delta_5(\otr)\,\,.
\ee

\begin{Theorem}
\label{noswirl}
Let $u$ be a bounded weak solution of the Navier-Stokes equations in
$\rthmiz$. Assume that $u$ is axi-symmetric with no swirl.
Then $u(x,t)=(0,0,b_3(t))$ for some bounded measurable function
$b_3\colon(-\i,0)\to R$.
\end{Theorem}

\begin{myproof}
The idea of the proof is the same as in the two-dimensional case.
By the results of Section~\ref{nsreg}, we have $|\nabla_x^k u|\le C_k$
in $R^3\times (-\i,0)$, and this implies that $\otr$ is bounded
in $\rthmiz$. Let $M_1=\sup_{\rthmiz}(\otr)$ and assume $M_1>0$.
Applying
Lemma~\ref{weakharnack} to equation \eqref{ote5}, considered
as an equation in $R^5\times(-\i,0)$, we see that $\otr\ge M_1/2$
in arbitrarily large parabolic balls (with suitably chosen centers).
However, this would mean that $\q \om \t$ is unbounded, a contradiction.
Therefore $M_1\le 0$. In the same way we show that
$M_2=\inf_{\rthmiz}\otr\ge 0$, and hence $\q \om \t$ vanishes identically.
For axi-symmetric vector fields with no swirl this means
that $\om=0$ and the proof is finished by  again applying the
Liouville theorem to the system $\curl u= 0, \,\,\div u = 0\,$. 
\end{myproof}
The validity of Theorem~\ref{noswirl} in the absence of the ``no swirl"
 assumption is still an open problem. The following theorem, however, is a partial result
in that direction:
\begin{Theorem}
\label{axisym}
Let $u$ be a bounded weak solution of the Navier-Stokes equations
in $\rthmiz$. Assume that $u$ is axi-symmetric and, in addition,
satisfies 
\be
\label{decayassumption}
|u(x,t)|\le\frac{C}{\sqrt{x_1^2+x_2^2}}\quad\mbox{in $\rthmiz$}.
\ee
Then $u=0$ in $\rthmiz$.
\end{Theorem} 

\begin{myproof}
We will use the cylindrical coordinates $(r,\t,z)$ given by 
\eqref{cylindrical}. We set $f=\rut$ and recall that
\be
\label{fequation}
f_t+\q u r f_{,\,r}+\q u z f_{,\,r}=\Delta f - \frac2{r}f_{,\,r}\,\,.
\ee
For $\lambda>0$ we let $\fl(x,t)=f(\lambda x,\lambda^2 t)$ and 
$\ul(x,t)=\lambda u(\lambda x,\lambda^2 t)$.
We note that $\fl$ again satisfies \eqref{fequation} with $u$ replaced
by $\ul$, a consequence of the fact that $\ul$ satisfies Navier-Stokes.
%For $\d>0$ we let
%\be
%\label{odef}
%\od=\{x\in R^3,\,\,\sqrt{x_1^2+x_2^2}\ge\d\}\times(-\i,0)\,\,.
%\ee   
%An important point is that under our assumptions we have, 
%for each $\d>0$
Under our assumptions we have
\begin{eqnarray}
|\fl| & \le & C\quad \mbox{in $\rthmiz\,$ uniformly in $\lambda>0$, and} \label{fbound} \\
|\ul| & \le & \frac{C}{r} \quad \mbox{in $\rthmiz\,$ uniformly in $\lambda>0$.}\label{ubound}
\end{eqnarray}
Let $M=\sup_{\rthmiz}f$. We will show that $M\le 0$. Arguing by contradiction,
let us assume that $M>0$.
Let us fix some $\d>0$. (It is instructive to think of $\delta$ as being small, 
although one can also take $\d=1$, for example.)
By re-scaling $f\to\fl$ we can move points where $\fl$ is ``almost equal to $M$" 
close to the $x_3-$axis.
Using this and Lemma~\ref{weakharnack}, we see that for any (large) $T_1>0$, $L>0$ and $R>0$ 
and any (small) $\ve>0$ we can find $\lambda>0$ such that 
$\fl\ge M-\ve$ in a space-time region $\ro$ of the form
\be
\ro=\{x\in R^3, \d\le \r\le R, -L+\bxt\le x_3\le L+\bxt\}\times(\bar t-T_1,\bar t)\,\,.
\ee
Consider a smooth axi-symmetric cut-off function $\vf(x,t)$ supported in
\be
\label{support}
\{x\in R^3,\,\r\le R,-L+\bxt\le x_3\le L+\bxt\}\times(\bar t - T_1,\bar t)
\ee
such that $\vf=1$ in
\be
\label{areaofone}
\{x\in R^3,\,\r\le R-1,-L+1+\bxt\le x_3\le L-1+\bxt\}\times(\bar t-T_1+1,\bar t -1)
\ee
and, moreover, $|\vf_t|\le1$, $|\vf_{,\,r}|\le 1$ and $ |\vf_{,\,z}|\le 1$ everywhere.
(A natural choice is, for example, $\vf(r,z,t)=\xi(r)\eta(z)\zeta(t)$ for suitable
functions $\xi,\eta,\zeta$ of one variable.)
Multiplying the equation for $\fl$ by $\vf$ and integrating over space-time,
we obtain
\be
\label{integrals}
\itix{(\fl_t+\q\ul r \fl_{,\,r}+\q \ul z\fl_{,\,z}-\Delta \fl)\vf}=
\itix{-\frac2 r\fl_{,\,r}\vf}\,\,.
\ee
This  equality will be shown to be impossible when $M>0$.
 In the integral on the left-hand side  of 
\eqref{integrals} one can change $\fl$ to $\fl-M$ and integrate by parts
to obtain
\be
\label{leftside}
\itix{(\fl-M)(-\vf_t-\ul\nabla\vf-\Delta\vf)}=I+II+III\,\,.
\ee
We have $\fl-M=O(\ve)$ in $\ro$ and therefore, if we allow correction terms of size 
$O(\ve)$, we can  restrict the spatial integration in these integrals to the
region $\{\r\le\d\}$. Using \eqref{fbound} and \eqref{ubound}, it is not hard to see
that 
\be
\label{leftbound}
|I|\le CL\d^2+O(\ve),\quad |II|\le C\d T_1+O(\ve)\quad\mbox{and $|III|\le C\d^2T_1+O(\ve)$.}
\ee
(We remind the reader that the value of $C$ can change from one expression to another.)
On the other hand, the right-hand side of \eqref{integrals} can be written as follows:
\be
\label{leftside1}
\itix{-\frac2 r\fl_{,\,r}\vf}=4\pi\int_{-\i}^0\int_{-\i}^{\i}\int_0^\i\fl\vf_{,\,r}\,dr\,dz\,dt
\ee
The key point then is that $\fl$ vanishes at the $x_3-$axis and is equal to
$M+O(\ve)$ on most of the support of $\vf_{,r}$. It is easy to check that the last integral in \eqref{leftside1}
is equal to
\be
\label{leftside2}
-4\pi M\int_{-\i}^{\i}\int_{-\i}^{\i}\vf(0,0,x_3,t)\,dx_3\,dt + O(\ve) \le -8\pi MLT_1+CT_1+CL+O(\ve)\,.
\ee
For $M>0$, this leads to a contradiction to \eqref{leftbound} and \eqref{integrals} if $L$ and $T_1$ are sufficiently large and
$\ve$ is sufficiently small. 
We have proved that $\sup f \le 0$. It follows in a similar way  that $\inf f\ge 0$ and therefore $f$ must vanish.
This means that the solution $u$ is swirl-free and we can apply Theorem~\ref{noswirl}
to conclude that $u=0$ in $\rthmiz$.   
\end{myproof}

\section{Singularities and ancient solutions}\label{axisymsing}
We will now consider the consequences of an assumption that a singularity exists in a solution of the Cauchy problem 
for the Navier-Stokes equations \eqref{NSeq} and \eqref{NSic}. We aim to show that singularities
generate  bounded {\it ancient solutions}, which are solutions defined in $R^n\times(-\i,0)$. 
More precisely, an {\it  ancient weak solution} of the Navier-Stokes
equations is a weak solution defined in $\rnmiz$, and $u$ is an {\it  ancient mild solution}
if there is a sequence $T_l\to-\i$ such that $u(\,\cdot\,,T_l)$ is well-defined and $u$ is a mild solution
of the Cauchy problem in $R^n\times(T_l,0)$ with initial datum $u(\,\cdot\,,T_l)$. %zastavka Sept 13, 1pm
(We remark that even if $u$ is a bounded weak solution of Navier-Stokes in $\rnmiz$, the function
$u(\,\cdot\,,t)$ may not be well-defined for each $t$, see Section~\ref{nsreg}. On the other hand 
$u(\,\cdot\,,t)$ is well defined for almost every $t$ for any $u\in\lxt\i\i(\rnmiz)$.)
\begin{Lemma}
\label{limitofmild}
Assume that $u_l$, $l=1,2,\dots$ is a sequence of bounded mild solutions of Navier-Stokes defined
in $R^n\times(T_l,0)$ (for some initial data) with a uniform bound $|u_l|\le C$, and $T_l\searrow-\i$.
Then we can choose a subsequence such that along the subsequence the $u_l$ converge locally uniformly
in $\rnmiz$ to an ancient mild solution $u$ satisfying $|u|\le C$ in $\rnmiz$.
\end{Lemma}

\begin{myproof}
This is an easy consequence of the results in Section~\ref{nsreg}.
\end{myproof}

\begin{Remark}
\label{usefulremark}
Another easy result, which is nevertheless a useful addendum to the Liouville theorems of Section~\ref{liouville}
is the following:
A bounded ancient mild solution $u(x,t)$ of the Navier-Stokes equations which is of the form $u(x,t)=b(t)$ 
is constant (independent of $t$).
\end{Remark}
We leave the proof of the last statement to the reader as a simple exercise.

Recall from Section~\ref{nsreg} that for any $u_0\in \lirn$ the Cauchy problem \eqref{NSeq}, \eqref{NSic}
has a unique local-in-time mild solution $u$.  
Assume now that the mild solution develops a singularity in finite time, and that $(0,T)$ is its  maximal time 
interval of existence. Let $h(t)=\sup_{x\in R^n}|u(x,t)|$. 
By a classical result of Leray (\cite{Leray}) we have 
 \be
 \label{growth}
 h(t)\ge\frac{\ve_1}{\sqrt{T-t}}
 \ee
 for some $\ve_1>0$. Let $H(t)=\sup_{0\le s\le t}h(s)$. It is easy to see that there exists a sequence
 $t_k\nearrow T$ such that $h(t_k)=H(t_k)$. Let us choose a sequence of numbers $\g_k\searrow 1$.
 Let  $N_k=H(t_k)$ and choose $x_k\in R^n$ such that 
 $M_k=|u(x_k,t_k)|\ge N_k/\g_k$.
 Let us set
 \be
 \label{scaling}
 \vk(y,s)=\frac{1} {M_k} u(x_k+\frac {y}{M_k},t_k + \frac{s}{M_k^2})\,.
 \ee
 The functions $\vk$ are defined in $R^n\times(A_k,B_k)$, with
 $A_k=-M_k^2 t_k$ and $B_k=M_k^2(T-t_k)\ge\ve_1^2\g_k^2$, and satisfy
 \be
 \label{vkbound}
 \mbox{$|\vk|\le \g_k$ in $R^n\times(A_k,0)$ and $|\vk(0,0)|=1$}.
 \ee
 Also, $\vk$ are mild solutions of the Navier-Stokes equations in
 $R^n\times(A_k,0)$ with initial data $v_0^{(k)}(y)=\frac1 {M_k}u_0(x_k+\frac {y}{M_k})$. 
 By Lemma~\ref{limitofmild}, there is a subsequence of $\vk$ converging to an
 ancient mild solution $v$ of the Navier-Stokes equations. By our construction, we have
 $|v|\le 1$ in $R^n\times(-\i,0)$ and $|v(0,0)|=1$.
 
 We have proved the following statement:
 \begin{Proposition}
 \label{principle}
 A finite-time singularity arising from a mild solution generates a bounded ancient mild solution which is not 
 identically zero.
 \end{Proposition}
 Without further information about the situation at hand,  the proposition may not be very useful. 
 By itself, the existence of non-zero bounded ancient solutions is not surprising. (Consider constants, for example.)
 However, if (non-zero) constant solutions can be excluded (for example by a scale-invariant estimate) and a Liouville-type
 theorem for ancient solutions is available, then finite-time singularities can be ruled out. 
 
A simple example of such a situation arises in the context of the Ladyzhen{-}skaya-Prodi-Serrin regularity criterion.
Assume that a finite $T>0$ is the maximal time of existence of a mild solution (with a suitable initial condition).
Let $p,q\in(1,\i)$ with $n/p+2/q=1$. Then $||u||_{\lxts p q(R^n\times(0,T))}=+\i$. To see this, it is enough to
note that if the $\lxts p q -$norm of $u$ was finite, the function $v$ constructed by the above procedure
would have to vanish identically a.\ e.\, due to the invariance of the $\lxts p q-$norm under the scaling used in the procedure,
along with the fact that the finiteness of the $\lxts p q-$norm implies its ``local smallness". But $v$ has to be smooth
(by the results of Section~\ref{nsreg}) and $|v(0,0)|=1$, a contradiction.

A more interesting application of the procedure gives Theorems~\ref{typetwospace} and \ref{typetwotime} below, 
which can be thought of as
generalizations of recent results in \cite{ChenStrainTsaiYau}.

\begin{Theorem}
\label{typetwospace}
Let $u$ be an axi-symmetric vector field in $R^3\times(0,T)$ which 
belongs to $\lxt\i\i(R^3\times(0,T'))$ for each $T'<T$.
Assume that $u$ is a weak solution of the Navier-Stokes equations
in $R^3\times(0,T)$ and that 
\be
\label{assumption}
\mbox{$|u(x,t)|\le \frac{C}{\r}$ in $R^3\times(0,T)$.}
\ee
Then $|u|\le M=M(C)$ in $R^3\times(0,T)$. Moreover, $u$ is a mild solution
of the Navier-Stokes equations (for a suitable initial datum).
\end{Theorem}
\begin{Remark}
By the results of Section~\ref{nsreg} regarding mild solutions we see that
$u$ is in fact smooth in $\rtht$ with pointwise bounds on all derivatives
in $R^3\times(\tau,T)$ for any fixed $\tau>0$.
\end{Remark}

\begin{myproof}
We first prove the statement assuming that $u$ is a mild solution (for a suitable
initial datum). This situation is in fact the main point of the theorem.
The fact that we can weaken the assumptions from mild solutions to weak solutions
 in the formulation
of the theorem (while keeping the other assumptions the same) is only of 
marginal interest.

Arguing by contradiction, let us assume that $u$ is a mild solution which
is bounded in $\rthtp$ for each $T'<T$ and develops a singularity at time $T$.
We now use the re-scaling procedure  described in the paragraph preceding
Proposition~\ref{principle} to %zastavka2
 construct a  bounded ancient mild solution
$v$. Let $x_k$ and $M_k$ be as in the construction. We will
write $x_k=(x'_k,x_{3k})$, with $x'_k=(x_{1k},x_{2k})$.
An obvious consequence of assumption~\ref{assumption} is that
$|x'_k|\le \frac C {M_k}$. This implies that the functions $\vk(y,s)$
are axi-symmetric with respect to an axis parallel to the
$y_3-$axis and at distance at most $C$ from it. Therefore we can
assume (by passing to a suitable subsequence first) that the limit function 
$v$ is axi-symmetric with respect to a suitable
axis. Moreover, since assumption~\eqref{assumption} is scale-invariant,
it will  again be satisfied (in suitable coordinates) by $v$. Applying Theorem~\ref{axisym}
and using \eqref{assumption} we see that $v=0$. On the other hand,
$|v(0,0)|=1$, a contradiction.
This finishes the main part of the proof. 

It remains to show that, under
the assumptions of the theorem, $u$ is a mild solution. To do this we inspect
the decomposition of $u$ constructed in Lemma~\ref{decomposition} with
$f_k=-u_k u$. Using the decay of the kernel \eqref{decay} and of the heat kernel,
 it is easy to check that, under the assumption \eqref{assumption}, all the
terms in the decomposition $u=v+w+b$ will again satisfy \eqref{assumption}.
It follows easily that $b$ must vanish and therefore $u$ is a mild solution.  
\end{myproof}

Theorem~\ref{typetwospace} can be used to prove the following result:
\begin{Theorem}
\label{typetwotime}
Let $u$ be an axi-symmetric vector field in $R^3\times(0,T)$ which 
belongs to $\lxt\i\i(R^3\times(0,T'))$ for each $T'<T$.
Assume that $u$ is a weak solution of the Navier-Stokes equations
in $R^3\times(0,T)$ satisfying
\be
\label{assumption1}
\mbox{$|u|\le \frac{C}{\sqrt{T-t}}$  in $R^3\times(0,T)$.}
\ee
In addition, assume that there exists some $R_0>0$ such that
\be
\label{assumption2}
\mbox{for $\r\ge R_0$ and $0<t<T$ we have $|u(x,t)|\le {\frac{C}{\r}}\,\,$,}
\ee
as is for example the case when $u$ is a mild solution with initial datum $u_0$
decaying sufficiently fast at $\i$.

Then $|u|\le M=M(C)$ in $R^3\times(0,T)$. Moreover, $u$ is a mild solution
of the Navier-Stokes equations (for a suitable initial datum).
\end{Theorem}
We remark that the statement  fails, for trivial reasons,  if we drop 
assumption \eqref{assumption2}. (Consider   $u(x,t)=b(t)$.)
The fact that \eqref{assumption2} is satisfied when $u_0$ decays sufficiently fast
at $\i$ (e.\ g.\ when it is compactly supported) follows for example 
from \cite{Brandolese1, Brandolese2}. 

\begin{myproof}
We have seen in the proof of Theorem~\ref{typetwospace} that \eqref{assumption2}
implies that $u$ is a mild solution for a suitable initial datum and is therefore  smooth
in open subsets of $\rtht$. 
We define
\be
\label{fquant}
f(x,t)=|x'| \,\,|u(x,t)|=\r\,\,\,|u(x,t)|,
\ee
where, as above, $x'=(x_1,x_2)$. By Theorem~\ref{typetwospace}, it is enough to prove
that $f$ is bounded in $\rtht$. 
 Let $h(t)=\sup_{\rth} f(x,t)$, $H(t)=\sup_{0\le\tau\le t}h(\tau)$. Assume $f$ is not bounded
  and choose $t_k\nearrow T$ and $x_k\in\rth$ such that
  $M_k=f(x_k,t_k)=h(t_k)=H(t_k)\nearrow \infty$. 
  Let $\lambda_k=|x'_k|$ and, for $y\in\rth,\,\,s\in(-T\lambda_k^{-2},0)$, define 
  \be
  \label{vkdef}
  \vk(y,s)=\vk(y',y_3,s)=\lambda_k u(\lambda_k y', \lambda_k y_3+x_{3k},T+\lambda_k^2 s)\,\,.
  \ee  
We note that the sequence ${\lambda_k}$ is bounded due to \eqref{assumption2}.
Set $s_k=-(T-t_k)\lambda_k^{-2}$.
Since \eqref{assumption1} is invariant under the Navier-Stokes scaling,
the functions $\vk$ satisfy
\be
\label{assumptionvk}
|\vk|\le \frac C{\sqrt{-s}}\quad\mbox{in $\rth\times(-T\lambda_k^{-2},0)$\,,}
\ee
where $C$ is the same as in \eqref{assumption1}\,\,.

Moreover, from the construction we have
\be
\label{vkbound1}
|\vk(y,s)|\le \frac{M_k}{|y'|} \quad\mbox{in $\rth\times(-T\lambda_k^{-2},s_k).$}
\ee
Note also that by the elementary inequality $\min(1/a,1/b)\le2/(a+b)\,\,,$ estimates
 \eqref{assumptionvk} and \eqref{vkbound1}
 imply
 \be
 \label{vkbound2}
 |\vk(y,s)|\le\frac{2CM_k}{M_k\sqrt{-s}+C|y'|}\quad\mbox{in $\rth\times(-T\lambda_k^{-2},s_k).$}
 \ee
 Let 
$\gamma\subset\rth$ be the unit circle $\{y\in\rth,\,|y'|=1,y_3=0\}\,$. 
We have, by construction, $|\vk(\,\cdot\,,s_k)|\big|_\gamma=M_k$ which, together with 
\eqref{assumptionvk} shows that $s_k\ge -C^2M_k^{-2}$.

Therefore, roughly speaking, as $k\to\infty$, the sequence $\vk$ blows up 
along $\gamma$. If we knew that the $\vk$ satisfied local energy estimates with bounds independent of $k$,
the blow-up along $\gamma$ would be in contradiction with the partial
regularity theory in \cite{CKN}, since the one-dimensional
Hausdorff measure of the blow-up set must be zero. One can in fact work along these lines
and finish the proof, but the procedure is not simple.

One can alternatively finish the proof by another scaling argument (one could do both scalings in one step,
but the two-step procedure seems to be more transparent): Denoting by $e_1$ the vector
$(1,0,0)$,  for $x\in\rth$ and $\tau\in(A_k,0]$ where $A_k=M_k^2(-T\lambda_k^{-2}-s_k)$, we define
\be
\label{wkdef}
\wk(x,\tau)=\frac1 {M_k}\vk(e_1+\frac x {M_k}, s_k+\frac \tau {M_k^2})\,.
\ee
We will consider the cylinders
\be
\label{cylinder}
\ck=\{x\in\rth,\,\, \sqrt{(x_1+M_k)^2+x_2^2}\le\frac{M_k} 2\}\,.
\ee
It follows from our definitions that 
\be
\label{wkbound}
\mbox{$|\wk(0,0)|=1$ and $|\wk(x,\tau)|\le 2$ in $(\rth\setminus\ck)\times(A_k,0)\,.$}
\ee 
Note also that \eqref{vkbound2} implies
\be
\label{wkbound2}
|\wk(x,\tau)|\le\frac{2CM_k}{M_k\sqrt{-\tau}+C\sqrt{(x_1+M_k)^2+x_2^2}}\quad\mbox{in $\ck\times(A_k,0)$}
\ee
and that ~\eqref{assumptionvk} implies 
\be
\label{wkbound3}
|\wk(x,\tau)|\le\frac C {\sqrt{-\tau}}\quad \mbox{in $\rth\times(A_k,0)$}.
\ee
Since the functions $\wk$ are mild solutions of the Navier-Stokes equations
in $(A_k,0)$ (for suitable re-scalings of the initial datum $u_0$), in view of bound \eqref{wkbound3}
we can choose a subsequence of the sequence $\wk$, which we again denote by $\wk$,
such that the $\wk$ converge uniformly on compact subsets of $\rth\times(-\i,0)$ to an ancient
mild solution $w$. In view of \eqref{wkbound} we have $|w|\le 2$ in
$\rth\times(-\i,0)$. Moreover, since the solutions $\vk$ are axi-symmetric and $M_k\nearrow\i$,
it is easy to see that $w$ is independent of the $x_2-$variable.
Applying Theorem~\ref{nstd} and Remark~\ref{usefulremark} to the field $(w_1,w_3)$, we conclude
that $(w_1,w_3)$ must vanish identically, and this easily implies that $w=0$ in $\rth\times(-\infty,0)$.
This would give a contradiction with $|\wk(0,0)|=1$ if we could prove that
$\wk(0,0)\to w(0,0)$, which is not immediately obvious since our bound of $\sup_x|\wk(x,\tau)|$ may not be
uniform as $\tau\to 0$. However, by \eqref{wkbound} the only possible problem may occur  
due to the contribution from the cylinder $\ck$. In the cylinder we can use the bound \eqref{wkbound2} 
to show that the contribution of the dangerous part of $\wk$ to the representation
formula \eqref{rep} is negligible (in the limit $k\to\infty$).
Applying the representation formula \eqref{rep} in $\rth\times(-1,0)$ with $\wk(x,-1)$ as initial
datum and $f_{jl}=-\wk_l \wk_j$ and using the bound \eqref{wkbound2} together with the decay of the kernel 
\eqref{decay}, one sees that it is enough to estimate the integral
\be
\label{last}
I(M)=\int_{-1}^0\int_{-\i}^\i\int_{|x'|\le \frac M 2}\,\,\,\frac1{(\sqrt{-\tau}+\frac{|x'|}{M})^2}\,\,\,\frac 1{(\frac {M^2} 4 + x_3^2)^2}\,\,\,dx' \,dx_3\,d\tau\,\,.
\ee 
An easy calculation shows that $I(M)\to0$ as $M\to\i$. This shows that the contribution from the region
where $|\wk|\ge2$ to the
representation formula \eqref{rep} (with $f_{jl}=-\wk_l\wk_j$) is negligible (in the limit $k\to\infty$) and therefore (by \eqref{estimates1})
the sequence $\wk$ converges to $w$ uniformly in $\bar B(0,1)\times[-1,0]$. 
Therefore $|w(0,0)|=1$, which gives the sought-after contradiction.

\end{myproof}

\end{document}